\documentclass{amsart}

\usepackage{amsmath,amsthm,amsfonts,amscd,amssymb}

\newtheorem{thm}{Theorem}[section]

\newtheorem{lem}[thm]{Lemma}

\theoremstyle{definition}

\theoremstyle{remark}
\newtheorem{rem}{Remark}[section]

\begin{document}

\title{Small zeros of hermitian forms over a quaternion algebra}
\author{Wai Kiu Chan and Lenny Fukshansky}

\address{Department of Mathematics and Computer Science, Wesleyan University, Middletown, CT 06459}
\email{wkchan@wesleyan.edu}
\address{Department of Mathematics, 850 Columbia Avenue, Claremont McKenna College, Claremont, CA 91711}
\email{lenny@cmc.edu}

\begin{abstract}
Let $D$ be a positive definite quaternion algebra over a totally real number field $K$, $F(\boldsymbol X,\boldsymbol Y)$ a hermitian form in $2N$ variables over $D$, and $Z$ a right $D$-vector space which is isotropic with respect to $F$. We prove the existence of a small-height basis for $Z$ over $D$, such that $F(\boldsymbol X,\boldsymbol X)$ vanishes at each of the basis vectors. This constitutes a non-commutative analogue of  a theorem of Vaaler \cite{vaaler:smallzeros2}, and presents an extension of the classical theorem of Cassels \cite{cassels:small} on small zeros of rational quadratic forms to the context of quaternion algebras.
\end{abstract}

\maketitle

\def\A{{\mathcal A}}
\def\AA{{\mathfrak A}}
\def\B{{\mathcal B}}
\def\C{{\mathcal C}}
\def\D{{\mathcal D}}
\def\F{{\mathcal F}}
\def\x{{\mathcal H}}
\def\I{{\mathcal I}}
\def\J{{\mathcal J}}
\def\K{{\mathcal K}}
\def\kk{{\mathfrak K}}
\def\L{{\mathcal L}}
\def\M{{\mathcal M}}
\def\mm{{\mathfrak m}}
\def\MM{{\mathfrak M}}
\def\NN{{\mathfrak N}}
\def\OO{{\mathfrak O}}
\def\O{{\mathcal O}}
\def\R{{\mathcal R}}
\def\s{{\mathcal S}}
\def\V{{\mathcal V}}
\def\UU{{\mathfrak U}}
\def\X{{\mathcal X}}
\def\Y{{\mathcal Y}}
\def\Z{{\mathcal Z}}
\def\H{{\mathcal H}}
\def\cee{{\mathbb C}}
\def\pee{{\mathbb P}}
\def\que{{\mathbb Q}}
\def\real{{\mathbb R}}
\def\zed{{\mathbb Z}}
\def\aaa{{\mathbb A}}
\def\hyp{{\mathbb H}}
\def\Bb{{\mathbb B}}
\def\Ff{{\mathbb F}}
\def\Nn{{\mathbb N}}
\def\kk{{\mathfrak K}}
\def\qbar{{\overline{\mathbb Q}}}
\def\kbar{{\overline{K}}}
\def\xbar{{\overline{x}}}
\def\ybar{{\overline{Y}}}
\def\kkbar{{\overline{\mathfrak K}}}
\def\ubar{{\overline{U}}}
\def\abar{{\overline{a}}}
\def\eps{{\varepsilon}}
\def\ahat{{\hat \alpha}}
\def\bhat{{\hat \beta}}
\def\gt{{\tilde \gamma}}
\def\h{{\tfrac12}}
\def\ba{{\boldsymbol a}}
\def\be{{\boldsymbol e}}
\def\bei{{\boldsymbol e_i}}
\def\bc{{\boldsymbol c}}
\def\bm{{\boldsymbol m}}
\def\bk{{\boldsymbol k}}
\def\bi{{\boldsymbol i}}
\def\bl{{\boldsymbol l}}
\def\bq{{\boldsymbol q}}
\def\bu{{\boldsymbol u}}
\def\bt{{\boldsymbol t}}
\def\bs{{\boldsymbol s}}
\def\bv{{\boldsymbol v}}
\def\bw{{\boldsymbol w}}
\def\bx{{\boldsymbol x}}
\def\bX{{\boldsymbol X}}
\def\bz{{\boldsymbol z}}
\def\bwy{{\boldsymbol y}}
\def\bY{{\boldsymbol Y}}
\def\bL{{\boldsymbol L}}
\def\bet{{\boldsymbol\eta}}
\def\bxi{{\boldsymbol\xi}}
\def\bo{{\boldsymbol 0}}
\def\bol{{\boldkey 1}_L}
\def\ep{\varepsilon}
\def\p{\boldsymbol\varphi}
\def\q{\boldsymbol\psi}
\def\Hf{H_{\fin}^{\O}}
\def\Hfo{H_{\fin}^{\O_1}}
\def\Hft{H_{\fin}^{\O_2}}
\def\Hfd{H_{\fin}^{O_D}}
\def\rank{\operatorname{rank}}
\def\aut{\operatorname{Aut}}
\def\lcm{\operatorname{lcm}}
\def\sgn{\operatorname{sgn}}
\def\spn{\operatorname{span}}
\def\md{\operatorname{mod}}
\def\Norm{\operatorname{Norm}}
\def\dim{\operatorname{dim}}
\def\det{\operatorname{det}}
\def\Vol{\operatorname{Vol}}
\def\rk{\operatorname{rk}}
\def\ord{\operatorname{ord}}
\def\ker{\operatorname{ker}}
\def\div{\operatorname{div}}
\def\Gal{\operatorname{Gal}}
\def\Tr{\operatorname{Tr}}
\def\nn{\operatorname{N}}
\def\inf{\operatorname{inf}}
\def\fin{\operatorname{fin}}
\def\Gr{\operatorname{Gr}}
\def\Mat{\operatorname{Mat}}
\def\GL{\operatorname{GL}}

\section{Introduction}
\label{intro}

Let $F(X_1,\dots,X_N)$ be a quadratic form in $N \geq 2$ variables with rational coefficients, and suppose that $F$ is isotropic over $\que$. In his celebrated 1955 paper \cite{cassels:small}, J. W. S. Cassels proved that in this case there must exist a non-trivial rational zero of $F$ of small height, where the explicit bound on the height is $\ll_N H(F)^{(N-1)/2}$; here $H(F)$ stands for the height of the quadratic form $F$, to be defined below. In an addendum to the same paper Cassels demonstrated an example due to M. Kneser, which shows that the exponent $(N-1)/2$ on $H(F)$ cannot be improved. Cassels' result has since been generalized and extended by a variety of authors to many different contexts. Notably, S. Raghavan \cite{raghavan} has proved a similar result for quadratic and hermitian forms with coefficients over a number field, A. Prestel \cite{prestel} extended Cassels' theorem to rational function fields, and A. Pfister \cite{pfister} proved it over algebraic function fields (see for instance \cite{me:quad} for additional more recent bibliography).

It appears, however, that the case of hermitian forms has not been studied further since Raghavan's paper. In this note we consider the situation of an isotropic hermitian form on a vector space over a positive definite quaternion algebra, defined over a totally real number field, and prove that there exists a basis for this space consisting entirely of small-height zeros of our hermitian form. More precisely, our main result is as follows, where height  functions $h$, $H_{\inf}$, and $H^{\O}$, the height on $D$ with respect to an order $\O$, are to be defined in section~\ref{heights} below.

\begin{thm} \label{main} Let $D = \binom{\alpha,\beta}{K}$ be a positive definite quaternion algebra over a totally real number field $K$, where $\alpha, \beta$ are totally negative algebraic integers in $K$. Let $\O$ be an order in $D$. Let $N \geq 2$ be an integer, and let $Z \subseteq D^N$ be an $L$-dimensional right $D$-subspace, $1 \leq L \leq N$. Let $F(\bX,\bY) \in D[\bX,\bY]$ be a hermitian form in $2N$ variables, and assume that $F$ is isotropic on $Z$. Then there exists a basis $\bwy_1,\dots,\bwy_L$ for $Z$ over $D$ such that $F(\bwy_n) := F(\bwy_n,\bwy_n) = 0$ for all $1 \leq n \leq L$ and
\begin{equation}
\label{mn1}
h(\bwy_1) \leq \A_K(N, L,\alpha,\beta) H_{\inf}(F)^{\frac{4L-1}{2}} H^{\O}(Z)^4,
\end{equation}
and
\begin{equation}
\label{mn2}
h(\bwy_1) h(\bwy_n) \leq \A_K(N, L,\alpha,\beta)^2 H_{\inf}(F)^{4L-1} H^{\O}(Z)^8,
\end{equation}
where the constant $\A_K(N, L,\alpha,\beta)$ is defined in (\ref{main_const}) below.
\end{thm}

All notation used in Theorem \ref{main} is defined in section~\ref{heights}. The investigation of small-height linearly independent zeros of quadratic forms goes back to the works of H. Davenport \cite{davenport}, J. H. H. Chalk \cite{chalk}, R. Schulze-Pillot \cite{schulze}, H. P. Schlickewei \cite{schlickewei}, H. P. Schlickewei and W. M. Schmidt \cite{schmidt:schlickewei}, and J. D. Vaaler \cite{vaaler:smallzeros}, \cite{vaaler:smallzeros2}, among others. Our Theorem \ref{main} is essentially a non-commutative analogue of a similar result by J. D. Vaaler in \cite{vaaler:smallzeros2} over a number field. We use Vaaler's result as one of the main tools in our argument. We use the height function machinery over quaternion algebras as defined by C. Liebend{\"o}rfer in \cite{liebendorf:1}. 

In section \ref{heights} we define our notation, introduce the necessary quadratic and hermitian forms, and develop heights in both commutative and non-commutative settings. The proof of our main result, which we present in section~\ref{proof}, uses comparison inequalities between heights over a quaternion algebra and heights over its ground number field, which we prove in section~\ref{height_comp}. We believe that these comparison lemmas will be useful for future work in non-commutative Diophantine analysis with height functions. We discuss the optimality of our bounds in Remark~\ref{opt_rem} at the end of section~\ref{proof}.
\bigskip

\section{Heights and quadratic forms}
\label{heights}

We start with some notation. Let $K$ be a number field of degree $d$ over $\que$, $O_K$ its ring of integers, $M(K)$ its set of places, $\Delta_K$ its discriminant, and let us write $\Nn$ for the norm from $K$ to $\que$. For each place $v \in M(K)$ we write $K_v$ for the completion of $K$ at $v$ and let $d_v = [K_v:\que_v]$ be the local degree of $K$ at $v$, so that for each $u \in M(\que)$
\begin{equation}
\sum_{v \in M(K), v|u} d_v = d.
\end{equation}

\noindent
For each place $v \in M(K)$ we define the absolute value $|\ |_v$ to be the unique absolute value on $K_v$ that extends either the usual absolute value on $\real$ or $\cee$ if $v | \infty$, or the usual $p$-adic absolute value on $\que_p$ if $v|p$, where $p$ is a rational prime. Then for each non-zero $a \in K$ the {\it product formula} reads
\begin{equation}
\label{product_formula}
\prod_{v \in M(K)} |a|^{d_v}_v = 1.
\end{equation} 

\noindent
We extend absolute values to vectors by defining the local heights. Let $N \geq 1$, and for each $v \in M(K)$ define a local height $H_v$ on $K_v^N$ by
$$H_v(\bx) = \max_{1 \leq i \leq N} |x_i|_v,$$
and for each $v | \infty$ define another local height $\H_v$ on $K_v^N$ by
$$\H_v(\bx) = \left( \sum_{i=1}^N |x_i|_v^2 \right)^{1/2}.$$
for each $\bx \in K_v^N$. Then we define two global height function on $K^N$:
\begin{equation}
H(\bx) = \prod_{v \in M(K)} H_v(\bx)^{d_v/d},\ \H(\bx) = \prod_{v \nmid \infty} H_v(\bx)^{d_v/d} \times \prod_{v | \infty} \H_v(\bx)^{d_v/d}
\end{equation}
for each $\bx \in K^N$. Notice that due to the normalizing exponent $1/d$, our global height functions are absolute, i.e. for points over $\qbar$ their values do not depend on the field of definition. This means that if $\bx \in \qbar^N$ then $H(\bx)$ and $\H(\bx)$ can be evaluated over any number field containing the coordinates of $\bx$.

We also define an {\it inhomogeneous} height function on vectors by
 \begin{equation}
h(\bx) = H(1,\bx),
\end{equation}
hence $h(\bx) \geq H(\bx)$ for each $\bx \in \qbar^N$. In fact, the values of $H$ and $h$ are also related in the following sense: for each $\bx \in K^N$, there exists $a \in K$ such that $a\bx \in O_K^N$ and
\begin{equation}
\label{H_to_h}
H(\bx) = h(a \bx).
\end{equation}

\smallskip

We will also define two different height functions on matrices. First, let $B$ be an $N \times N$ matrix with entries in $K$, then we can view $B$ as a vector in $K^{N^2}$ and write $H(B)$ to denote the height of this vector. In particular, if $B$ is a symmetric matrix, then
$$Q(\bX,\bY) = \bX^t B \bY$$
is a symmetric bilinear form in $2N$ variables over $K$, and 
$$Q(\bX) := Q(\bX,\bX) =  \bX^t B \bX$$
is the associated quadratic form in $N$ variables. We define $H(Q)$, the height of such quadratic and bilinear forms, to be $H(B)$.
\smallskip

The second height we define on matrices is the same as height function on subspaces of $K^N$. Let $X = (\bx_1 \dots \bx_L)$ be an $N \times L$ matrix of rank $L$ over $K$, $1 \leq L \leq N$. Define
\begin{equation}
\label{matrix_ht}
\H(X) = \H(\bx_1 \wedge \dots \wedge \bx_L).
\end{equation}
For each $v | \infty$, the Cauchy-Binet formula guarantees that
\begin{equation}
\label{cauchy_binet}
\H_v(X) = |\det (X^* X)|_v^{1/2},
\end{equation}
where $X^*$ is the complex conjugate transpose of $X$. On the other hand, $\bx_1 \wedge \dots \wedge \bx_L$ can be identified with the vector $\Gr(X)$ of {\it Grassmann coordinates} of $X$ under the canonical embedding into $K^{\binom{N}{L}}$. Namely, let $\I$ be the collection of all subsets $I$ of $\{1,...,N\}$ of cardinality $L$, then $|\I| = \binom{N}{L}$. For each $I \in \I$, write $X_I$ for the $L \times L$ submatrix of $X$ consisting of all those rows of $X$ which are indexed by $I$. Define
\begin{equation}
\label{GR}
\Gr(X) = (\det (X_I))_{I \in \I} \in K^{\binom{N}{L}}.
\end{equation}
By our remark above, $\H(X) =  \H(\Gr(X))$. Now let $V \subseteq K^N$ be an $L$-dimensional subspace, $1 \leq L \leq N$. Choose a basis $\bx_1,...,\bx_L$ for $V$ over $K$, and let $X = (\bx_1\ ...\ \bx_L)$ be the corresponding $N \times L$ basis matrix. Define height of $V$ to be
$$H(V) := \H(X).$$
This height is well defined, since it does not depend on the choice of the basis for $V$: let $\bwy_1,...,\bwy_L$ be another basis for $V$ over $K$ and $Y = (\bwy_ 1 \dots \bwy_L)$ the corresponding $N \times L$ basis matrix, then there exists $C \in \GL_L(K)$ such that $Y = XC$, and so
$$\bwy_1 \wedge \dots \wedge \bwy_L = (\det C)\ \bx_1 \wedge \dots \wedge \bx_L,$$
hence, by the product formula $\H(\bwy_1 \wedge \dots \wedge \bwy_L) = \H(\bx_1 \wedge \dots \wedge \bx_L)$.

It will be convenient for us to define certain field constants that we use in our inequalities. Following \cite{vaaler:smallzeros}, first define for every $v \in M(K)$
\begin{equation}
\label{r_v}
r_v(L) = \left\{ \begin{array}{ll}
\pi^{-1/2} \Gamma \left( \frac{L}{2}+1 \right)^{1/L} & \mbox{if $v$ is real} \\
(2\pi)^{-1/2} \Gamma \left( L+1 \right)^{1/2L} & \mbox{if $v$ is complex} \\
1 & \mbox{if $v \nmid \infty$},
\end{array}
\right.
\end{equation}
and let
\begin{equation}
\label{CKL}
C_K(L) = 2 |\Delta_K|^{1/2d} \prod_{v \in M(K)} r_v(L)^{d_v/d}, 
\end{equation}
for each $L \geq 1$. Finally let
\begin{equation}
\label{BKL}
B_K(L) = 2^{L+1} C_K(1)^2 C_K(L-1)^{2(L-1)}.
\end{equation}
\bigskip

We can also extend the height machinery to the context of quaternion algebras, using the approach of \cite{liebendorf:1}. Let $K$ as above be a totally real number field, then $K$ has precisely $d$ archimedean places $v_1,\dots,v_d$, corresponding to the embeddings 
\begin{equation}
\label{embed}
\sigma_n : K \rightarrow K_{v_n} \cong \real
\end{equation}
for all $1 \leq n \leq d$. Then for each $a \in K$, $|a|_{v_n} = |\sigma_n(a)|$, where $|\ |$ stands for the usual absolute value on $\real$. We will also write $a^{(n)}$ for the algebraic conjugate $\sigma_n(a)$ of $a \in K$ under $\sigma_n$. Let $\alpha,\beta \in O_K$ be totally negative, meaning that $\alpha^{(n)} := \sigma_n(\alpha) < 0$ and  $\beta^{(n)} := \sigma_n(\beta) < 0$ for all $1 \leq n \leq d$. Let $D = \binom{\alpha,\beta}{K}$ be a positive definite quaternion algebra over~$K$, generated by the elements $i,j,k$ which satisfy the following relations:
\begin{equation}
\label{quat_rel}
i^2=\alpha,\ j^2=\beta,\ ij=-ji=k,\ k^2=-\alpha\beta.
\end{equation}
As a vector space, $D$ has dimension four over $K$, and $1,i,j,k$ is a basis. From now on we will fix this basis, and thus will always write each element $x \in D$ as
$$x = x(0)+x(1)i+x(2)j+x(3)k,$$
where $x(0),x(1),x(2),x(3) \in K$ are respective components of $x$, and the standard involution on $D$ is conjugation: 
$$\xbar= x(0)-x(1)i-x(2)j-x(3)k.$$
We define trace and norm on $D$ by
$$\Tr(x) = x+\xbar = 2x(0),\ \nn(x) = x\xbar = x(0)^2 - \alpha x(1)^2 - \beta x(2)^2 + \alpha \beta x(3)^2.$$
The algebra $D$ is said to be positive definite because the norm $\nn(x)$ is given by a positive definite quadratic form. In fact, since the norm form $\nn(x)$ is positive definite, $D_{v_n} := D \otimes_K K_{v_n}$ is isomorphic to the real quaternion $\hyp=\real + \real i + \real j + \real k$ for each $1 \leq n \leq d$. Hence each embedding $\sigma_n$ of $K$, $1 \leq n \leq d$, induces an embedding $\sigma_n : D \to D_{v_n}$, given by 
$$\sigma_n(x) = x(0)^{(n)} + x(1)^{(n)}i + x(2)^{(n)}j + x(3)^{(n)}k.$$
From now on we will write $x^{(n)}$ for $\sigma_n(x)$. Then the local norm at each archimedean place is also a positive definite quadratic form over the respective real completion~$K_{v_n}$:
\begin{eqnarray*}
& & \nn^{(n)}(x) = x^{(n)}\xbar^{(n)} \\
& & = \left(x(0)^{(n)}\right)^2 - \alpha^{(n)} \left(x(1)^{(n)}\right)^2 - \beta^{(n)} \left(x(2)^{(n)}\right)^2 + \alpha^{(n)} \beta^{(n)} \left(x(3)^{(n)}\right)^2,
\end{eqnarray*}
for each $1 \leq n \leq d$. We now have archimedean absolute values on $D$, corresponding to the infinite places $v_1,\dots,v_d$ of $K$: for each $x \in D$, define 
$$|x|_{v_n} = \sqrt{\nn^{(n)}(x)},$$
for every $1 \leq n \leq d$. It will be convenient to define
\begin{equation}
\label{s_ab_v}
s_{v_n}(\alpha,\beta) = \max \{1,|\alpha|_{v_n},|\beta|_{v_n},|\alpha \beta|_{v_n}\}^{\frac{1}{2}},\ t_{v_n}(\alpha,\beta) = \min \{1,|\alpha|_{v_n},|\beta|_{v_n},|\alpha \beta|_{v_n}\}^{\frac{1}{2}},
\end{equation}
for each $1 \leq n \leq d$, and also let
\begin{equation}
\label{s_ab}
s(\alpha,\beta) = \prod_{n=1}^d s_{v_n}(\alpha,\beta),\ t(\alpha,\beta) = \prod_{n=1}^d t_{v_n}(\alpha,\beta).
\end{equation}
Since local norm forms are positive definite, we immediately have the following inequalities:
\begin{equation}
\label{loc_ineq}
t_{v_n}(\alpha,\beta) \max_{0 \leq m \leq 3} |x(m)|_{v_n} \leq\ |x|_{v_n} \leq 2 s_{v_n}(\alpha,\beta) \max_{0 \leq m \leq 3} |x(m)|_{v_n}.
\end{equation}
Now, generalizing notation of \cite{liebendorf:1}, we can define an infinite homogeneous height on $D^N$ by
\begin{equation}
\label{Hinf}
H_{\inf}(\bx) = \left( \prod_{n=1}^d \max_{1 \leq l \leq N} |x_l|_{v_n} \right)^{1/d},
\end{equation}
and define an infinite inhomogeneous height on $D^N$ by
\begin{equation}
\label{hinf}
h_{\inf}(\bx) = H_{\inf}(1,\bx),
\end{equation}
for every $\bx \in D^N$. Clearly, $H_{\inf}(\bx) \leq h_{\inf}(\bx)$. The infinite height takes into account the contributions at the archimedean places. As in \cite{liebendorf:1}, we also define its counterpart, the finite height. Let us once and for all fix an order $\O$ in $D$; our definition will be with respect to the order $\O$, and this height will be denoted by $\Hf$. Specifically, for each $\bx \in \O^N$, let
\begin{equation}
\label{Hfin_O}
\Hf(\bx) = [\O : \O x_1 + \dots + \O x_N]^{-1/4d}.
\end{equation}
This is well defined, since $\O x_1 + \dots + \O x_N$ is a left submodule of $\O$. Now we can define the global homogeneous height on $\O^N$ by
\begin{equation}
\label{HD}
H^{\O}(\bx) = H_{\inf}(\bx) \Hf(\bx),
\end{equation}
and the global inhomogeneous height by
\begin{equation}
\label{hD}
h(\bx) := H_{\inf}(1,\bx) \Hf(1,\bx) = h_{\inf}(\bx) \geq H^{\O}(\bx),
\end{equation}
since $\O + \O x_1 + \dots + \O x_N = \O$. To extend this definition to $D^N$, notice that for each $\bx \in D^N$ there exists $a \in O_K$ such that $a\bx \in \O^N$, and define $H^{\O}(\bx)$ to be $H^{\O}(a \bx)$ for any such~$a$. This is well defined by the product formula, and $H^{\O}(\bx t) = H^{\O}(\bx)$ for all $t \in D^{\times}$.
\smallskip

We will now define height on the set of proper right $D$-subspaces of $D^N$, again following \cite{liebendorf:1}. Recall that $D$ splits over $E=K(\sqrt{\alpha})$, meaning that there exists a $K$-algebra homomorphism $\rho: D \to \Mat_{22}(E)$, given~by
\begin{equation}
\label{rho}
\rho(x(0)+x(1)i+x(2)j+x(3)k) = \left( \begin{matrix} x(0)+x(1)\sqrt{\alpha} & x(2)+x(3)\sqrt{\alpha} \\ \beta(x(2)-x(3)\sqrt{\alpha}) & x(0)-x(1)\sqrt{\alpha} \end{matrix} \right),
\end{equation}
so that $\rho(D)$ spans $\Mat_{22}(E)$ as an $E$-vector space (see Proposition 13.2a (p. 238) and Exercise 1 (p. 240) of \cite{pierce}). This map extends naturally to matrices over $D$. Let $Z \subseteq D^N$ be an $L$-dimensional right vector subspace of $D^N$, $1 \leq L < N$. Then there exists an $(N-L) \times N$ matrix $C$ over $D$ with left row rank $N-L$ such that $Z$ is the solution space of the linear system $C\bX = \bo$.  Define
\begin{equation}
\label{Hinf_C}
H_{\inf}(C) = \left( \prod_{n=1}^d \left| \det\left( \rho(CC^*) \right) \right|_{v_n} \right)^{1/4d},
\end{equation}
where $C^*$ is the conjugate transpose of $C$. The analogue of Cauchy-Binet formula works here as well (see (2.7) and (2.8) of \cite{liebendorf:1}, as well as Corollary 1 of \cite{liebendorf:2}), and so we have an alternative formula:
\begin{equation}
\label{Hinf_C1}
H_{\inf}(C) = \left( \prod_{n=1}^d \sum_{C_0} \left| \det\left( \rho(C_0) \right) \right|^2_{v_n} \right)^{1/2d},
\end{equation}
where the sum is taken over all $(N-L) \times (N-L)$ minors $C_0$ of $C$. Also define
\begin{equation}
\label{Hfin_C}
\Hf(C) = [\O^{N-L} : C(\O^N)]^{-1/4d},
\end{equation}
where $C$ is viewed as a linear map $\O^N \to \O^{N-L}$. Then we can define
\begin{equation}
\label{H_Z}
H^{\O}(Z) = H^{\O}(C) := H_{\inf}(C) \Hf(C).
\end{equation}
This definition does not depend on the specific choice of such matrix~$C$. By the duality principle proved in \cite{liebendorf:3},
\begin{equation}
\label{H_Z_dual}
H^{\O}(Z) = H^{\O}(Z^{\perp}),
\end{equation}
where $Z^{\perp} = \{ \bwy \in D^N : \bx^* \bwy = 0\ \forall\ \bx \in Z \}$. This means that if $\bx_1,\dots,\bx_L$ is a basis for $Z$ over $D$ and $X = (\bx_1 \dots \bx_L)$ is the corresponding basis matrix, then 
\begin{equation}
\label{H_Z_dual1}
H^{\O}(Z) = H^{\O}(X) := \left( [\O^L : X^t(\O^N)]^{-1} \prod_{n=1}^d \left| \det\left( \rho(X^*X) \right) \right|_{v_n} \right)^{1/4d},
\end{equation}
completely analogous to the definition of the height $H^{\O}(C)$ in (\ref{H_Z}); here $X^t$ is viewed as a linear map $\O^N \to \O^{N-L}$.

It will also be convenient to define a map $[\ ]: D \to K^4$, given by 
$$[x] = (x(0),x(1),x(2),x(3)),$$
for each $x = x(0)+x(1)i+x(2)j+x(3)k \in D$. This map obviously extends to $[\ ]: D^N \to K^{4N}$, given by $[\bx] = ([x_1],\dots,[x_N])$ for each $\bx = (x_1,\dots,x_N) \in D^N$. Clearly this is a bijection; in fact, it is an isomorphism of $K$-vector spaces, and we will write $[\ ]^{-1}$ for its inverse. 
\smallskip

By analogy with heights over $D$, we will also write
$$H_{\inf}(\bx) = \prod_{v | \infty} H_v(\bx)^{d_v/d},\ H_{\fin}(\bx) = \prod_{v \nmid \infty} H_v(\bx)^{d_v/d},$$
for every $\bx \in K^N$. Then by Lemma 2.1 of \cite{liebendorf:1}, for every $\bx \in O_K^N$ we have
\begin{equation}
\label{hqf1}
H_{\fin}(\bx) = \left[ O_K : O_K x_1 + \dots + O_K x_N \right]^{-1/d}.
\end{equation}
Also, if $V$ is an $L$-dimensional subspace of $K^N$ and $C$ is any $(N-L) \times N$ matrix over $O_K$ of rank $1 \leq L < N$, viewed as a linear map $O_K^N \to O_K^{N-L}$, such that $V = \{ \bx \in K^N : C\bx = \bo \}$, let us write
$$H_{\inf}(C) = \prod_{v | \infty} \H_v(C)^{d_v/d},\ H_{\fin}(C) = \prod_{v \nmid \infty} H_v(C)^{d_v/d},$$
and then by Lemma 2.1 and Proposition 2.4 of \cite{liebendorf:1}, we have
\begin{equation}
\label{hqf2}
H_{\fin}(C) = \left[ O_K^{N-L} : C(O_K ^N) \right]^{-1/d}.
\end{equation}
This means that the definitions over $K$ and over $D$ are really analogous.


Now let $F(\bX,\bY) \in D[\bX,\bY]$ be a hermitian form in $2N$ variables with coefficients in $D$, so that $F(a \bx,\bwy) = \bar{a} F(\bx,\bwy)$ and $F(\bwy,\bx) = \overline{F(\bx,\bwy)}$ for each $a \in D$ and $\bx,\bwy \in D^N$. We also write $F(\bX)$ for $F(\bX,\bX)$, then $F(\bx) \in K$ for any $\bx \in D^N$. Let us also write $\Ff = (f_{ml})$ for the $N \times N$ coefficient matrix of $F$, then $f_{ml} = \overline{f_{lm}}$ for each $1 \leq l,m \leq N$, and $F(\bX,\bY) = \bX^t \Ff \bY$. Same way as for quadratic and bilinear forms over $K$, we will talk about the height of the hermitian form $F$ over $D$, where by $H^{\O}(F)$ (respectively,  $H_{\inf}(F)$, $\Hf(F)$) we will always mean $H^{\O}(\Ff)$ (respectively,  $H_{\inf}(\Ff)$, $\Hf(\Ff)$), viewing $\Ff$ as a vector in $D^{N^2}$. We define the corresponding bilinear form $B$ over $K$ by taking the trace of $F$, i.e. $B([\bX],[\bY]) = \Tr(F(\bX,\bY))$, then the associated quadratic form
\begin{equation}
\label{trace_form}
Q([\bX]) := B([\bX],[\bX])
\end{equation}
in $4N$ variables over $K$ is equal to $2F(\bX)$. Therefore $F(\bx) = 0$ for some $\bx \in D^N$ if and only if $Q([\bx]) = 0$. Write $\Bb$ for the $4N \times 4N$ symmetric matrix of $B$ over $K$, then each entry of $\Ff$ corresponds to a $4 \times 4$ block in $\Bb$. Specifically, if $f_{ml} = f_{ml}(0)+ f_{ml}(1)i+ f_{ml}(2)j+ f_{ml}(3)k \in D$, then the corresponding block in $\Bb$ is of the form
\begin{equation}
\label{trace_matrix}
\Bb(f_{ml}) := \left( \begin{matrix}
2f_{ml}(0)&2\alpha f_{ml}(1)&2\beta f_{ml}(2)&-2\alpha \beta f_{ml}(3)\\
-2\alpha f_{ml}(1)&-2\alpha f_{ml}(0)&-2\alpha \beta f_{ml}(3)&2\alpha \beta f_{ml}(2)\\
-2\beta f_{ml}(2)&2\alpha \beta f_{ml}(3)&-2\beta f_{ml}(0)&-2\alpha \beta f_{ml}(1)\\
2\alpha \beta f_{ml}(3)&-2\alpha \beta f_{ml}(2)&2\alpha \beta f_{ml}(1)&2\alpha \beta f_{ml}(0)
\end{matrix} \right),
\end{equation}
so $\Bb = (\Bb(f_{ml}))_{1 \leq m,l \leq N}$, and $Q(\bz) = \bz^t \Bb \bz$ for each $\bz \in K^{4N}$. As defined before, we will write $H(Q)$ (respectively,  $H_{\inf}(Q)$, $H_{\fin}(Q)$) for $H(\Bb)$ (respectively,  $H_{\inf}(\Bb)$, $H_{\fin}(\Bb)$), viewed as a vector in $K^{16N^2}$.

Finally, we define the constant that appears in the upper bounds of our main result, Theorem \ref{main}. 
Let $\Delta_{\O}$ be the discriminant of the order $\O$, which is the ideal in $O_K$ generated by all the elements of the form
$$\det \left( \Tr (\omega_h \omega_n) \right)_{0 \leq h,n \leq 3} \in O_K,$$
where $\omega_0,\dots,\omega_3$ are in $\O$, and let
\begin{equation}
\label{O_const}
\MM(\O) := \max \left\{ \frac{\Nn(\Delta_{\O})^{1/2}}{\Nn(4 \alpha \beta)}, \frac{\Nn(4 \alpha \beta)}{\Nn(\Delta_{\O})^{1/2}} \right\},
\end{equation}
where $\Nn$ stands for the norm from $K$ to $\que$. Then define
\begin{equation}
\label{main_const}
\A_K(N, L,\alpha,\beta) = 2^{\frac{20 L - 3}{2}} N^{4L-1} \sqrt{B_K(4L)}\ \MM(\O)^{4(N-L)} \left( \frac{s(\alpha,\beta)^{4L}}{t(\alpha,\beta)^{\frac{4L-1}{2}}} \right),
\end{equation}
where the field constant $B_K(L)$ is defined in (\ref{BKL}), $\MM(\O)$ is as in (\ref{O_const}), and $s(\alpha,\beta)$, $t(\alpha,\beta)$ are defined in (\ref{s_ab}). We are now ready to proceed.
\bigskip

\section{Height comparison lemmas}
\label{height_comp}

In this section we will derive some inequalities between heights over~$K$ and over~$D$, which we later use to prove our main result. We start with the following simple lemma.

\begin{lem} \label{ht_comp1} For each $\bx \in \O^N$,
\begin{equation}
\label{comp1}
t(\alpha,\beta) H([\bx]) \leq H_{\inf}(\bx) \leq h_{\inf}(\bx) = h(\bx) \leq 2 s(\alpha,\beta) h([\bx]).
\end{equation}
\end{lem}

\proof
For each $v|\infty$, (\ref{loc_ineq}) implies that
$$t_v(\alpha,\beta) H_v([\bx]) \leq \max_{1 \leq l \leq N} |x_l|_v \leq 2 s_v(\alpha,\beta) H_v([\bx]),$$
and $h_{\inf}(\bx) = h(\bx)$ by (\ref{hD}).
\endproof

Similarly, we can obtain inequalities between height of a hermitian form over $D$ and its associated trace form.

\begin{lem} \label{ht_comp2} Let $F$ be a hermitian form over $D$ and let $Q$ be its associated trace form over $K$, as in (\ref{trace_form}). Then
\begin{equation}
\label{comp2}
\frac{t(\alpha,\beta)}{2 s(\alpha,\beta)^2} H(Q) \leq H_{\inf}(F),\ H^{\O}(F) \leq 2^{\frac{d+1}{d}} s(\alpha,\beta) \Nn(\alpha \beta)^{\frac{1}{d}} \NN(\O) H(Q),
\end{equation}
where $\Nn$ stands for the norm on $K$, and
\begin{equation}
\label{NN}
\NN(\O) = \min \left\{ |\Nn(\gamma)|^{\frac{1}{d}} : \gamma \in O_K \text{ is such that } \gamma i, \gamma j, \gamma k \in \O \right\}.
\end{equation}
\end{lem}

\proof
As in section~\ref{heights} above, let $\Ff = (f_{ml})$ be the $N \times N$ coefficient matrix of $F$. There exists $a \in O_K$ such that simultaneously all the entries of the matrix $a \Ff$ are in $\O$ and all the entries of the matrix $a \Bb$ of the corresponding bilinear trace form of $aF$ as defined in section~\ref{heights} above are in~$O_K$. Notice that $H^{\O}(a \Ff) = H^{\O}(\Ff)$ and $H(a \Bb) = H(\Bb)$ by the product formula. Hence we may assume without loss of generality that all entries of $\Ff$ are in $\O$ and all entries of $\Bb$ are in~$O_K$. Now for each $1 \leq m,l \leq N$ and $v \in M(K)$ such that $v | \infty$,
\begin{eqnarray*}
H_v(\Bb(f_{ml})) & \leq & 2 s_v(\alpha,\beta)^2 \max_{0 \leq k \leq 3} |f_{ml}(k)|_v \\ 
& \leq & \frac{2s_v(\alpha,\beta)^2}{t_v(\alpha,\beta)} |f_{ml}|_v \leq \frac{4s_v(\alpha,\beta)^3}{t_v(\alpha,\beta)} H_v(\Bb(f_{ml})),
\end{eqnarray*}
by (\ref{loc_ineq}), and so
$$\frac{t_v(\alpha,\beta)}{2 s_v(\alpha,\beta)^2} H_v(\Bb) \leq \max_{1 \leq m,l \leq N} |f_{ml}|_v \leq 2 s_v(\alpha,\beta) H_v(\Bb),$$
meaning that
\begin{equation}
\label{c1}
\frac{t(\alpha,\beta)}{2 s(\alpha,\beta)^2} H_{\inf}(Q) \leq H_{\inf}(F) \leq 2 s(\alpha,\beta) H_{\inf}(Q).
\end{equation}
There exists $\gamma \in O_K$ such that $\gamma i, \gamma j, \gamma k \in \O$. Out of all such elements $\gamma$ let us pick the one with the minimal norm $|\Nn(\gamma)|$. Notice that
$$H_{\inf}(\gamma F) = |\Nn(\gamma)|^{\frac{1}{d}} H_{\inf}(F) = \NN(\O) H_{\inf}(F),$$
where $\NN(\O)$ is as in (\ref{NN}), and so (\ref{c1}) becomes
\begin{equation}
\label{c1.1}
\frac{t(\alpha,\beta)}{2 s(\alpha,\beta)^2} H_{\inf}(Q) \leq H_{\inf}(F)  = \NN(\O)^{-1} H_{\inf}(\gamma F) \leq 2 s(\alpha,\beta) H_{\inf}(Q).
\end{equation}
\smallskip

Now we need to obtain a similar inequality for the finite heights. For each $1 \leq m,l \leq N$, write $\Bb(f_{ml})_{n,h}$ for the $nh$-th entry of the $4 \times 4$ matrix $\Bb(f_{ml})$, $0 \leq n,h \leq 3$. Then, by (\ref{hqf1}) we have
\begin{equation}
\label{c2}
H_{\fin}(\Bb) = \left[ O_K : \sum_{m,l=1}^N \sum_{n,h=0}^3 O_K \Bb(f_{ml})_{n,h} \right]^{-1/d} \leq 1.
\end{equation}

On the other hand, 
\begin{eqnarray}
\label{c3}
\Hf(\gamma F) & = & \left[ \O : \sum_{m,l=1}^N \O \gamma f_{ml} \right]^{-\frac{1}{4d}}  \nonumber \\
& \leq & \left[ \O : \sum_{m,l=1}^N \O \gamma f_{ml}(0) + \O \gamma i f_{ml}(1) + \O \gamma j f_{ml}(2) + \O \gamma k f_{ml}(3)  \right]^{-\frac{1}{4d}} \nonumber \\
& \leq & \left[ \O : \sum_{m,l=1}^N \sum_{h=0}^3 \O f_{ml}(h) \right]^{-\frac{1}{4d}} = \left[ O_K : \sum_{m,l=1}^N \sum_{h=0}^3 O_K f_{ml}(h) \right]^{-\frac{1}{d}} \nonumber \\
& \leq & \left( 2 \Nn(\alpha \beta \right))^{\frac{1}{d}} \left[ O_K : \sum_{m,l=1}^N \sum_{n,h=0}^3 O_K \Bb(f_{ml})_{n,h} \right]^{-\frac{1}{d}} \nonumber \\
& = & \left( 2 \Nn(\alpha \beta \right))^{\frac{1}{d}} H_{\fin}(\Bb),
\end{eqnarray}
where the identity in the third line of (\ref{c3}) follows by (2.23) of \cite{liebendorf:1}. Now observe that $H^{\O}(F) = H^{\O}(\gamma F)$ by the product formula, and so (\ref{comp2}) follows by combining (\ref{c1.1}), (\ref{c2}), and (\ref{c3}).
\endproof
\smallskip

\begin{rem} \label{42_counter_ex} Notice that (\ref{c2}) cannot be replaced with an inequality of the form $H_{\fin}(\Bb) \ll \Hf(F)$, which means that $H_{\inf}(F)$ in the upper bound of the first inequality of (\ref{comp2}) cannot be replaced with $H^{\O}(F)$. Consider the following example. Let $D = \binom{-1,-1}{\que}$, and let $\O$ be the order spanned by $\frac{1+i+j+k}{2},i,j,k$ over $\zed$. Let $n \in \zed$, and let $F$ be the hermitian form over $D$ given by the matrix
$$\left( \begin{matrix} 0 & i+nj \\ -i-nj & 0 \end{matrix} \right).$$
Then
$$\Hf(F) = \left[ \O : \O(i+nj) \right]^{-1/4} = \left( \Nn(i+nj)^2 \right)^{-1/4} = (1+n^2)^{-1/2}.$$
On the other hand, $\Bb(f_{11})$ is the $4 \times 4$ zero matrix, while
$$\Bb(f_{12}) = \left( \begin{matrix}
0 & -2 & -2n & 0 \\
2 & 0 & 0 & 2n \\
2n & 0 & 0 & -2 \\
0 & -2n & 2 & 0
\end{matrix} \right),$$
and so 
$$H_{\fin}(\Bb) = \left[ \zed : 2\zed + 2n\zed \right]^{-1} = 1/2.$$
Since the integer $n$ can be arbitrarily large, it is clear that we cannot have $H_{\fin}(\Bb) \ll \Hf(F)$ in this case.
\end{rem}
\smallskip


In the next lemma we will compare the heights of a $D$-subspace of $D^N$ with respect to two different orders, $\O_1$ and $\O_2$ in $D$.

\begin{lem} \label{ht_O_OD}  Let $\O_1$ and $\O_2$ be two orders in $D$, and let $Z \subseteq D^N$ be an $L$-dimensional right vector $D$-subspace of $D^N$, $1 \leq L \leq N$. Then
\begin{equation}
\label{O_OD}
 \M^{-(N-L)} H^{\O_1}(Z)  \leq H^{\O_2}(Z) \leq  \M^{N-L} H^{\O_1}(Z),
\end{equation}
where
\begin{equation}
\label{M_const}
\M = \M(\O_1,\O_2) := \max \left\{ \Nn \left( \Delta_{\O_1} \Delta_{\O_2}^{-1} \right)^{\frac{1}{2}}, \Nn \left( \Delta_{\O_2} \Delta_{\O_1}^{-1} \right)^{\frac{1}{2}} \right\}.
\end{equation}
\end{lem}

\proof
Let $C$ be an $(N-L) \times N$ matrix over $D$ with left row rank $N-L$ such that $Z$ is the solution space of the linear system $C\bX=\bo$, then we can view $C$ as a linear map $D^N \to D^{N-L}$ and
\begin{equation}
\label{ho12.1}
\Hfo(C) = [\O_1^{N-L} : C(\O_1^N)]^{-1/4d},\ \Hft(C) = [\O_2^{N-L} : C(\O_2^N)]^{-1/4d}.
\end{equation}
Let $\UU_1$ be the order ideal of the pair $\O_1, \O_2$, which is the product of the invariant factors of $\mathcal O_2$ in $\mathcal O_1$ (see for instance p. 49 of \cite{reiner}).  Then $\UU_1$ is a fractional ideal of $O_K$ generated by all the elements of the form $\det(\phi)$ as $\phi$ runs through all $K$-linear maps sending $\O_1$ into $\O_2$. In an analogous manner, let $\UU_2$ be the order ideal of the pair $\O_2, \O_1$. Then (\ref{ho12.1}) implies that
\begin{equation}
\label{ho12.2}
\Hfo(C) \leq  \Nn(\UU_1)^{N-L} \Hft(C),\ \Hft(C) \leq  \Nn(\UU_2)^{N-L} \Hfo(C).
\end{equation}
By definition of $\UU_1$ and $\UU_2$, it is easy to see that
$$\Delta_{\O_2} = \UU_1^2 \Delta_{\O_1},\ \Delta_{\O_1} = \UU_2^2 \Delta_{\O_2},$$
and so
\begin{equation}
\label{ho12.3}
 \Nn(\UU_1) = \Nn( \Delta_{\O_2} \Delta_{\O_1}^{-1} )^{\frac{1}{2}},\  \Nn(\UU_2) = \Nn( \Delta_{\O_1} \Delta_{\O_2}^{-1} )^{\frac{1}{2}}.
\end{equation}
Combining (\ref{ho12.2}) with (\ref{ho12.3}) finishes the proof.
\endproof
\smallskip


\begin{rem} \label{OD} Notice that if $\O_1 = \O_2$, then $\M=1$, and so there is equality throughout in~(\ref{O_OD}). In fact, given any two orders $\O_1$ and $\O_2$ in $D$, the existence of constants $C_1, C_2$ depending on $\O_1,\O_2,N,L$ such that
$$C_1 H^{\O_1}(Z) \leq H^{\O_2}(Z) \leq C_2 H^{\O_1}(Z)$$
for an $L$-dimensional right vector $D$-subspace $Z$ of $D^N$ was observed by Daniel Bertrand, and is discussed on p. 116 of \cite{liebendorf:1}.
\end{rem}
\smallskip

We will apply Lemma \ref{ht_O_OD} in the particular situation when $\O_1$ is the order $\O$ we picked in section~\ref{heights} above, and $\O_2$ is the order $O_D$ defined by
\begin{equation}
\label{eta}
O_D = \sum_{n=0}^3 O_K \eta_n, \text{ where } \eta_0 = 1,\ \eta_1 = i,\ \eta_2 = j,\ \eta_3 = k.
\end{equation}

\begin{lem} \label{ht_O_OD_1} Let $Z \subseteq D^N$ be an $L$-dimensional right vector $D$-subspace of $D^N$, $1 \leq L \leq N$. Then
\begin{equation}
\label{O_OD_1}
\MM(\O)^{-(N-L)} H^{\O}(Z)  \leq H^{O_D}(Z) \leq  \MM(\O)^{N-L} H^{\O}(Z),
\end{equation}
where $\MM(\O)$ is defined in (\ref{O_const}) above.
\end{lem}

\proof
A simple computation shows that $\Delta_{O_D}$ is generated by
$$\det \left( \Tr (\eta_h \eta_n) \right)_{0 \leq h,n \leq 3} = -16 \alpha^2 \beta^2,$$
and the lemma now follows from Lemma \ref{ht_O_OD}.
\endproof

Next we compare the height of a $D$-subspace of $D^N$ with respect to $O_D$ with the height of its image under~$[\ ]$.

\begin{lem} \label{ht_comp3}  Let $Z \subseteq D^N$ be an $L$-dimensional right vector $D$-subspace of $D^N$, $1 \leq L \leq N$, then $V_Z := [Z] \subseteq K^{4N}$ is a $4L$-dimensional $K$-subspace of $K^{4N}$, and
\begin{equation}
\label{comp3}
H(V_Z) = H^{O_D}(Z)^4,
\end{equation}
where $O_D$ is as in (\ref{eta}) above.
\end{lem}

\proof
Let $\sigma_1,\dots,\sigma_d$ be the $d$ real embeddings of $K$ as in (\ref{embed}), and define
\begin{equation}
\label{embed1}
\Sigma = (\sigma_1,\dots,\sigma_d): K \to \prod_{n=1}^d K_{v_n} \cong \real^d,
\end{equation}
which extends naturally to a map from $K^N$ to $\real^{Nd}$. Then the composition of $[\ ]$ with $\Sigma$ is an embedding of $D^N$ into $\real^{4Nd}$, and 
$$\Lambda_Z := \Sigma([Z \cap O_D^N]) = \Sigma(V_Z \cap O_K^{4N})$$
is a lattice of rank $4Ld$ in $\real^{4Nd}$, since $[O_D^N] = O_K^{4N}$. By a well-known theorem of W. M. Schmidt (see Theorem 1 on p. 435 of \cite{schmidt:67}), we can relate determinant of $\Lambda_Z$ to the height of $V_Z$. Specifically, as was worked out in (17) of \cite{me:number},
\begin{equation}
\label{schm1}
\det (\Lambda_Z) = |\Delta_K|^{L/2} H(V_Z)^d.
\end{equation}
On the other hand, C. Liebend\"{o}rfer in \cite{liebendorf:1} relates determinant of $\Lambda_Z$ to the height of $Z$. Lemma 3.2 of \cite{liebendorf:1} implies that
$$\det (\Lambda_Z) = \det(\Sigma[O_D])^L H(Z)^{4d},$$
where a straightforward extension of Lemma 3.1 of \cite{liebendorf:1} gives
$$\det(\Sigma[O_D]) = \sqrt{|\Delta_K|},$$
and so
\begin{equation}
\label{schm2}
\det (\Lambda_Z) = |\Delta_K|^{L/2} H(Z)^{4d}.
\end{equation}
It should be remarked that Lemmas 3.1 and 3.2 in \cite{liebendorf:1} are proved for the case when the number field $K$ is just $\que$, however the extensions of these lemmas to the case when $K$ is a totally real number field are essentially word for word: we are simply viewing $O_D$ as an $O_K$-module, which itself is a $\zed$-module. The exponent $d$ appearing on heights in our identities (\ref{schm1}) and (\ref{schm2}) is just a normalization due to the fact that our heights our absolute. Now combining (\ref{schm1}) with (\ref{schm2}) produces~(\ref{comp3}).
\endproof

\bigskip

\section{Main result}
\label{proof}

We are now ready to prove Theorem \ref{main}.

\proof[Proof of Theorem \ref{main}]
Let $Z \subseteq D^N$ be an $L$-dimensional right vector $D$-subspace of $D^N$, $1 \leq L \leq N$, and let $F$ be a hermitian form in $N$ variables over $D$ which is non-singular on $Z$. Then $V_Z = [Z] \subseteq K^{4N}$ is a $4L$-dimensional $K$-subspace of $K^{4N}$, and the corresponding quadratic form $Q$ is non-singular on $V_Z$. 

Notice that zeros of the hermitian form $F$ in $Z$ are in bijective correspondence with zeros of the quadratic form $Q$ in $V_Z$, i.e. $F(\bx)=0$ for some $\bx \in Z$ if and only if $Q([\bx])=0$. By formula (1.4) of \cite{vaaler:smallzeros2} (which is a corollary of Theorem 1 of \cite{vaaler:smallzeros2}) combined with our formula (\ref{H_to_h}), there exists a basis $\bx_1,\dots,\bx_{4L} \in O_K^{4N}$ for $V_Z$ such that $Q(\bx_l) = 0$ for all $1 \leq l \leq 4L$, $h(\bx_1) \leq h(\bx_l)$ and
\begin{equation}
\label{vaaler}
h(\bx_1) h(\bx_l) \leq B_K(4L) \left( 16 N^2 H(Q) \right)^{4L-1} H(V_Z)^2,
\end{equation}
for $1 \leq l \leq L$, where $B_K(L)$ is as in (\ref{BKL}). In particular,
\begin{equation}
\label{vaaler1}
h(\bx_1) \leq \sqrt{B_K(4L)} \left( 16 N^2 H(Q) \right)^{\frac{4L-1}{2}} H(V_Z),
\end{equation}
which is the analogue of Cassels' bound, and is sharp at least with respect to the exponent on $H(Q)$.

\begin{rem} \label{heights_rem} The height of the quadratic form $Q$ used by Vaaler in \cite{vaaler:smallzeros2} is slightly different from ours: he uses $L_2$-norms at the archimedean places instead of the sup-norms that we use for our height $H(Q)$. To compensate for this difference, we introduce the additional constant $16 N^2$ in front of $H(Q)$ in (\ref{vaaler}) and (\ref{vaaler1}). Our choice of the sup-norms at the archimedean places makes comparison inequalities of section~\ref{height_comp} (especially Lemma \ref{ht_comp2}) more natural and easier to work out.
\end{rem} 

Notice that $F([\bx_l]^{-1}) = 0$ for each $1 \leq l \leq 4L$, and there exist
$$1=l_1 < l_2 < \dots < l_L < 4L$$
such that $[\bx_{l_1}]^{-1},\dots,[\bx_{l_L}]^{-1}$ is a basis for $Z$ as a right $D$-vector space; we will write $\bwy_n = [\bx_{l_n}]^{-1}$ for each $1 \leq n \leq L$. Notice that in fact $\bwy_1,\dots,\bwy_L \in O_D^N$. Now we need to estimate the heights of these basis vectors, for which purposes we use the height comparison lemmas from section~\ref{height_comp}. Specifically, inequalities (\ref{mn1}) and (\ref{mn2}) follow from inequalities (\ref{vaaler1}) and (\ref{vaaler}), respectively, after an application of the height comparison inequalities presented in Lemmas \ref{ht_comp1} - \ref{ht_comp3} as follows: Lemma \ref{ht_comp1} produces an upper bound for each $h(\bwy_n)$ in terms of $h(\bx_{l_n})$; Lemma \ref{ht_comp2} is then used to bound $H(Q)$ from above in terms of $H_{\inf}(F)$; finally, using Lemma \ref{ht_comp3} we can express $H(V_Z)$ in terms of $H^{O_D}(Z)$, and then use Lemma \ref{ht_O_OD} to bound $H^{O_D}(Z)$ in terms of $H^{\O}(Z)$ for an arbitrary order $\O$. Hence combining (\ref{vaaler}) and (\ref{vaaler1}) with Lemmas \ref{ht_comp1} - \ref{ht_comp3}, we obtain inequalities (\ref{mn1}) and (\ref{mn2}). This completes the proof.
\endproof
\bigskip

\begin{rem} \label{basis_OD} Notice that the basis $\bwy_1,\dots,\bwy_L$ for $Z$ over $D$ we constructed in fact consists of vectors with all of their coordinates in the order~$O_D$. Moreover, for any order $\O'$ in $D$ there exists $b(\O') \in K$ such that the new basis vectors $b(\O') \bwy_1, \dots, b(\O') \bwy_L$ for $Z$ have all their coordinates in~$\O'$ while their projective heights $H^{\O}$ stay the same as those of $\bwy_1,\dots,\bwy_L$, respectively.
\end{rem}

\begin{rem} \label{opt_rem} It is not clear to us whether the upper bounds in Theorem \ref{main} are optimal or not. We should remark, however, that starting from the bounds of Theorem \ref{main} and applying our height comparison lemmas, one does retrieve a Cassels-type bound with correct exponents for the corresponding quadratic trace form over $K$ in $4N$ variables.  More precisely, assume for instance that the hermitian form $F$ has all of its coefficients in $O_D$, and at least one of these coefficients is equal to 1, so $H^{O_D}_{\fin}(F) = 1$ and thus $H_{\inf}(F) = H^{O_D}(F)$. Let $\bwy_1 \in Z$ be a non-trivial zero of $F$ guaranteed by (\ref{mn1}); let $\bx_1=[\bwy_1] \in V_Z$ be the corresponding zero of $Q$. Then combining (\ref{mn1}) with Lemmas \ref{ht_comp1} - \ref{ht_comp3}, we obtain
\begin{equation}
\label{opt}
H(\bx) \ll_{K,L,\alpha,\beta} H(Q)^{\frac{4L-1}{2}} H(V_Z).
\end{equation}
If we take $Z=D^N$, and so $V_Z=K^{4N}$, then $L=N$, $H(V_Z)=1$, and (\ref{opt}) becomes precisely a classical Cassels-type bound for $Q$ over $K$. Moreover, our exponent on $H_{\inf}(F)$ in (\ref{mn1}) is completely analogous to the bound obtained by Raghavan for hermitian forms over number fields (see Theorem 2 and remark at the bottom of page 114 in \cite{raghavan}).
\end{rem}
\bigskip

\bibliographystyle{plain}  
\bibliography{quaternion}        

\begin{thebibliography}{10}

\bibitem{cassels:small}
J.~W.~S. Cassels.
\newblock Bounds for the least solutions of homogeneous quadratic equations.
\newblock {\em Proc. Cambridge Philos. Soc.}, 51:262--264, 1955.

\bibitem{chalk}
J.~H.~H. Chalk.
\newblock Linearly independent zeros of quadratic forms over number fields.
\newblock {\em Monatsh. Math.}, 90(1):13--25, 1980.

\bibitem{davenport}
H.~Davenport.
\newblock Homogeneous quadratic equations. {P}repared for publication by {D}.
  {J}. {L}ewis.
\newblock {\em Mathematika}, 18:1--4, 1971.

\bibitem{me:number}
L.~Fukshansky.
\newblock Siegel's lemma with additional conditions.
\newblock {\em J. Number Theory}, 120(1):13--25, 2006.

\bibitem{me:quad}
L.~Fukshansky.
\newblock Small zeros of quadratic forms over $\overline{Q}$.
\newblock {\em Int. J. Number Theory}, 4(3):503--523, 2008.

\bibitem{liebendorf:1}
C.~Liebend{\"o}rfer.
\newblock Linear equations and heights over division algebras.
\newblock {\em J. Number Theory}, 105(1):101--133, 2004.

\bibitem{liebendorf:2}
C.~Liebend{\"o}rfer.
\newblock Heights and determinants over quaternion algebras.
\newblock {\em Comm. Algebra}, 33(10):3699--3717, 2005.

\bibitem{liebendorf:3}
C.~Liebend{\"o}rfer and G.~R{\'e}mond.
\newblock Duality of heights over quaternion algebras.
\newblock {\em Monatsh. Math.}, 145(1):61--72, 2005.

\bibitem{pfister}
A.~Pfister.
\newblock Small zeros of quadratic forms over algebraic function fields.
\newblock {\em Acta Arith.}, 79(3):221--238, 1997.

\bibitem{pierce}
R.~S. Pierce.
\newblock {\em Associative Algebras}.
\newblock Springer-Verlag, 1982.

\bibitem{prestel}
A.~Prestel.
\newblock On the size of zeros of quadratic forms over rational function
  fields.
\newblock {\em J. Reine Angew. Math.}, 378:101--112, 1987.

\bibitem{raghavan}
S.~Raghavan.
\newblock Bounds of minimal solutions of diophantine equations.
\newblock {\em Nachr. Akad. Wiss. Gottingen, Math. Phys. Kl.}, 9:109--114,
  1975.

\bibitem{reiner}
I.~Reiner.
\newblock {\em Maximal Orders}.
\newblock Academic Press, 1975.

\bibitem{schlickewei}
H.~P. Schlickewei.
\newblock Kleine nullstellen homogener quadratischer gleichungen.
\newblock {\em Monatsh. Math.}, 100(1):35--45, 1985.

\bibitem{schmidt:schlickewei}
H.~P. Schlickewei and W.~M. Schmidt.
\newblock Quadratic geometry of numbers.
\newblock {\em Trans. Amer. Math. Soc.}, 301(2):679--690, 1987.

\bibitem{schmidt:67}
W.~M. Schmidt.
\newblock On heights of algebraic subspaces and diophantine approximations.
\newblock {\em Ann. of Math. (2)}, 85:430--472, 1967.

\bibitem{schulze}
R.~Schulze-Pillot.
\newblock Small linearly independent zeros of quadratic forms.
\newblock {\em Monatsh. Math.}, 95(3):241--249, 1983.

\bibitem{vaaler:smallzeros}
J.~D. Vaaler.
\newblock Small zeros of quadratic forms over number fields.
\newblock {\em Trans. Amer. Math. Soc.}, 302(1):281--296, 1987.

\bibitem{vaaler:smallzeros2}
J.~D. Vaaler.
\newblock Small zeros of quadratic forms over number fields, {II}.
\newblock {\em Trans. Amer. Math. Soc.}, 313(2):671--686, 1989.

\end{thebibliography}

\end{document}